%% file: Dickman-distr.tex


 \magnification1050
 \def\oui{oui}
\ifx\textures\oui
 \input iten.tex
\input CrayolaColors 
 \long\def\rge#1{\Red#1\Black}

 \else
 \input itenn.tex

 \fi
\input option_keys.tex
\def\paradouze{oui}

\optionparag=1
\def\paradouze{oui}
\anglais

    \font\tenrsfs=rsfs7 at 10pt

    \font\sevenrsfs=rsfs7
    \font\sixrsfs=rsfs7 at 6pt
    
\newfam\rsfsfam\textfont\rsfsfam=\tenrsfs\scriptfont\rsfsfam=\sevenrsfs\scriptscriptfont\rsfsfam=\sixrsfs

\def\EE{{\bb E}}
\def\VV{{\bb V}}
 
\def\fl#1{\left\lfloor #1 \right\rfloor}

\ifx\montrerlabels\oui

\input montrerlabels.tex
\fi

\dimstand
\vsize235truemm
\hautspages{R. de la Bretche \& G. Tenenbaum}{Local limit theorems for a probabilistic model of the Dickman distribution}
\dateheure

\titrecentre{On strong and almost sure local limit theorems for a probabilistic model of the Dickman distribution}
\bigskip\medskip
\centerline{RŽgis de la Bretche \& GŽrald Tenenbaum}
\bigskip
{\leftskip75mm
\obeylines
\it To the memory of Jonas Kubilius,
 who stood on the bridge and invited us all.\par  }
\bigskip
{\eightpoint\leftskip1cm\rightskip1cm
\noi{\bf Abstract.}  \PAR
\medskip\noi
Let  $\{Z_k\}_{k\geqslant 1}$ denote  a sequence  of independent Bernoulli random variables defined by \hbox{$\PP(Z_k=1)=1/k=1-\PP(Z_k=0)$} $(k\geqslant 1)$ and  put $T_n:=\sum_{1\leqslant k\leqslant n}kZ_k$. It is then known that  $T_n/n$ converges weakly to a real random variable $D$ with density proportional to the Dickman function, defined by the delay-differential equation $u\varrho'(u)+\varrho(u-1)=0$ $(u>1)$ with initial condition $\varrho(u)=1$ $(0\leqslant u\leqslant 1)$. Improving on earlier work, we propose asymptotic formulae with remainders for the corresponding local and almost sure limit theorems, namely \smallskip
$$\sum_{m\geqslant 0}\abs{\PP(T_n=m)-{\e^{-\gamma}\over n}\varrho\Big({m\over n}\Big)}={2\log n\over \pi^2 n}\Big\{1+O\Big({1\over \log_2n}\Big)\Big\}\qquad (n\to\infty),$$ 
and
$$(\forall u>0) \sum_{\di{n\leqslant N}{T_n=\fl{un}}}1=\e^{-\gamma}\varrho(u)\log N+O\big((\log N)^{2/3+o(1)}\big)\quad\hbox{a.s.}\quad(N\to\infty),$$
where $\gamma$ denotes Euler's constant.
\medskip\noi
\medskip
\noi \bf 2020 Mathematics Subject Classification: \rm primary 60F15; secondary: 11N25, 60C05.\par }
\par 
\bigskip\bigskip
\paraun{Introduction and statement of results}
Dickman's function is defined on $[0,\infty[$ as the continuous solution to the delay-differential equation $u\varrho'(u)+\varrho(u-1)=0$ $(u>1)$ with initial condition $\varrho(u)=1$ $(0\leqslant u\leqslant 1)$. It is known (see, \eg, \citeplus{Te15}{th. III.5.10}) that $\int_0^\infty\varrho(u)\d u=\e^{\gamma}$, where $\gamma$ denotes Euler's constant. The Dickman distribution is defined as the law of a random variable $D$ on $[0,\infty[$ with density $$\varrho_0(u):=\e^{-\gamma}\varrho(u)Ê\quad (u\geqslant 0).$$
\par 
This law appears in a large variety of mathematical topics, such as (the following list being non limiting): \par\smallskip
{\leftskip5mm\rightskip5mm\hangafter1\hangindent8mm Number theory, in the context of friable integers,\note{That is integers free of large prime factors} after the seminal paper of Dickman \citer{Di30} : see \citer{Te15} for a expositary account; \par 
\hangafter1\hangindent8mm Random polynomials over finite fields : see, \eg,  Car~\citer{Ca87}, Manstavi\v cius \citer{Ma92}, Arratia, Barbour \& TavarŽ \citer{ABT93},  Knopfmacher \& Manstavi\v cius \citer{KM97};\par 
\hangafter1\hangindent8mm Random permutations: see in particular, Shepp \& Lloyd \citer{SL66}, Kingman \citer{Ki77}, Arratia, Barbour \& TavarŽ \citer{ABT00}, Manstavi\v cius \& Petuchovas \citer{MP17}.
\par 
In number theory, the Dickman function also appears in Billingsley's model \citer{Bi72} for the vector distribution of large prime factors of integers (see \citer{Te00} for an effective version) and in Kubilius' model\note{A probabilistic model for the uniform probability defined on the set of the first $N$ integers with $\sigma$-algebra comprising those events that can be defined by divisibility conditions involving solely small primes}: see Elliott \citer{El79-80}, Arratia, Barbour \& TavarŽ \citer{ABT99}, Tenenbaum \citer{Te99}, and \citeplus{Te15}{\S\thinspace III.6.5} for an expositary account.}
\par 
A simple probabilistic description of $D$ is provided by the almost surely convergent series $$\sum_{n\geqslant 1}\prod_{1\leqslant j\leqslant n}X_j,$$ where the $X_j$ are independent and uniform on $[0,1]$: see Goldie \& GrŸbel \citer{GG96}, Fill \& Huber~\citer{FH10}, Devroye \citer{De01}.
\par \goodbreak
There is a vast bibliography on the various probabilistic models of the Dickman distribution: see, \eg, Chen \& Hwang \citer{CH03}, Devroye \& Fawzi \citer{DF10},
Pinsky  \citer{Pi19}.
\par 
In 2002, Hwang \& Tsai \citer{HT02} used a simple model to show that, suitably normalized, the cost of Hoare's quickselect algorithm converges weakly to $D$. This  model may be described as follows: if $\{Z_k\}_{k\geqslant 1}$ denotes  a sequence  of independent Bernoulli random variables such that $\PP(Z_k=1)=1/k=1-\PP(Z_k=0)$ $(k\geqslant 1)$ and if $T_n:=\sum_{1\leqslant k\leqslant n}kZ_k$, then $T_n/n$ converges weakly to $D$, viz.
$$\lim_{n\to\infty}\PP(T_n\leqslant nu)=\e^{-\gamma}\int_0^u\varrho(v)\d v\qquad (u\geqslant 0).$$
\par 
 A strong local limit theorem was then obtained by Giuliano, Szewczak \&  Weber \citer{GSW18}, in the form 
$$v_n:=\sum_{m\geqslant 0}\abs{\PP(T_n=m)-{\e^{-\gamma}\over n}\varrho\Big({m\over n}\Big)}=o(1)\qquad (n\to\infty).\eqdef{vn}$$
\par  We propose a sharp estimate of the speed of convergence. Here and in the sequel, we let $\log_k$ denote the $k$-fold iterated logarithm.
\Propt{T1}{ We have
$${v_n={2\log n\over \pi^2 n}\Big\{1+O\Big({1\over \log_2n}\Big)\Big\}}\qquad (n\to\infty).\eqdef{encvn}$$}
\par 
 {This estimate may be put in perspective with the following result of Manstavi\v cius \citer{Ma92}. Let $\{X_k\}_{k\geqslant 1}$ denote  a sequence  of independent Poisson variables such that $\EE(X_k)=1/k$, and put $Y_n:=\sum_{1\leqslant k\leqslant n}kX_k$. Then \citeplus{Ma92}{cor. 2} readily yields the strong local limit theorem
$$\sum_{m\geqslant 0}\abs{\PP(Y_n=m)-{\e^{-\gamma}\over n}\varrho\Big({m\over n}\Big)}\ll{1\over n}\qquad (n\geqslant 1).\eqdef{fMa}$$
Thus, as may be expected, Poissonian approximations to the Bernoulli random variables $Z_k$ provide a closer model of the Dickman distribution. As a byproduct of \eqref{encvn} and \eqref{fMa}, we get an estimate of the total variation distance between $T_n$ and $Y_n$, viz.
$$\eqalign{d_{TV}(T_n,Y_n)&:=\sum_{m\geqslant 0}\abs{\PP(T_n=m)-\PP(Y_n=m)}=v_n+O\Big({1\over n}\Big)\cr&={ 2\log n\over \pi^2 n}\Big\{1+O\Big({1\over \log_2n}\Big)\Big\}.}$$
\par 
We also point, without details, to a recent estimate of Bhattacharjee \& Goldstein \citeplus{BG19}{th.1.1}, which provides a bound $\leqslant 3/(4n)$ for a smooth Wasserstein-type distance between $T_n/n$ and~$D$.}
\par \medskip
For $u>0$, let $\varepsilon_n$ denote a non-negative sequence tending to 0 at infinity, and let $\{m_n\}_{n\geqslant 1}$ denote a non-decreasing integer sequence such that $m_n=un+O(\varepsilon_nn)$ as $n\to\infty$. We may then define a sequence   of random variables $\{L_N(u)\}_{N=1}^{\infty}$ by the formula $$L_N(u):=\sum_{n\leqslant N,\,T_n=m_n}1.$$ By a complicated proof resting on a general correlation inequality, an almost sure local limit theorem is established in \citer{GSW18} assuming furthermore that $\{m_n\}_{n\geqslant 1}$ is strictly increasing: for any $u\geqslant 1$, the asymptotic formula  $L_N(u)\sim\e^{-\gamma}\varrho(u)\log N$ holds almost surely as $N\to\infty$.\note{The authors of \citer{GSW18} state  that this almost sure asymptotic formula holds for all $u>0$. However,  the requirement that $\{m_n\}_{n=1}^{\infty}$ should be strictly increasing is incompatible with the assumption $m_n\sim un$ if $u< 1$.} 
 \par 
The following result, proved by a simple, direct method, provides an effective version.
\Propt{T2}{Let $u\geqslant 1$, $\varepsilon_n=o(1)$ as $n\to\infty$, and let $\{m_n\}_{n\geqslant 1}$ denote a strictly increasing sequence of integers such that $m_n=un+O(\varepsilon_nn)$ $(n\geqslant 1)$. We have, almost surely,
$$L_N(u)=\bigg\{1+O\Big(\eta_N+{(\log_2N)^{1/2+o(1)}\over (\log N)^{1/3}}\bigg)\bigg\}\e^{-\gamma}\varrho(u)\log N,\eqdef{faLN}$$
where $\eta_N:=(1/\log N)\sum_{1\leqslant N}\varepsilon_n/n=o(1)$.
\par \smallskip
Furthermore, for any $u>0$, the formula $L_N(u)\sim\e^{-\gamma}\varrho(u)\log N$ holds almost surely  provided $$\vartheta_m:=|\{n\geqslant 1:m_n=m\}|=o(\log m)\qquad (m\to\infty),$$
and assuming only that the sequence $\{m_n\}_{n\geqslant 1}$ is  non-decreasing.  If $\vartheta_m\ll(\log m)^{\alpha}$ with $0\leqslant \alpha<1$,  the estimate \eqref{faLN} holds with  remainder $\ll\eta_N+1/(\log N)^{(1-\alpha)/3+o(1)}$.}
\medskip
We note that, for all $u>0$, the case $m_n:=\fl{un}$ is covered by the second part of the statement with~$\alpha=0$.\bigskip
\paraun{Proof of \ref{T1}}
Let $c$ be a large constant and put $M(n):=cn(\log n)/(\log_23n)$. We first show that the contribution to $v_n$ of those $m>M(n)$  is negligible. Indeed, since $\varrho(v)\ll v^{-v}$, we first have
$${1\over n}\sum_{m>M(n)}\varrho\Big({m\over n}\Big)\ll{1\over n}\sum_{m> M(n)}\e^{-m(\log_2n)/2n}\ll{1\over n}\cdot$$
Then, we have, for all $y\geqslant 0$
$$\sum_{m>  M(n) }\PP(T_n=m)\leqslant \e^{-yM(n)}\EE(\e^{yT_n})= \e^{-yM(n)} \prod_{1\leqslant k\leqslant n}\Big(1+{\e^{ky}-1\over k}\Big).\eqdef{majTngrand}$$
Selecting $y=(\log_2 n)/n$, we see that the last product is
$$\ll \exp\Big\{\int_0^n{\e^{yv}-1\over v}\d v\Big\}\ll n^{ 2/\log_2 n}, $$
hence the left-hand side of \eqref{majTngrand} is also $\ll1/n$, and we infer that
$$v_n=\sum_{1\leqslant m\leqslant M(n)}\abs{\PP(T_n=m)-{\e^{-\gamma}\over n}\varrho\Big({m\over n}\Big)}+O\Big({1\over n}\Big).\eqdef{est1vn}$$ \par 
\par 
 Recall the definition $\varrho_0(u):=\e^{-\gamma}\varrho(u)$ $(u\in\r)$ and let $I(s):=\int_0^1(\e^{vs}-1) \d v/v$ $(s\in\CC)$. 
From \citeplus{Te15}{th. III.5.10}, we know that
$$\wh{\varrho_0}(\tau):=\int_\r\e^{i\tau u}\varrho_0(u)\d u=\e^{I(i\tau)}\qquad (\tau\in\r).\eqdef{tfr}$$
Next, for $|\tau|<\pi$, we have
$$\eqalign{\EE\big(\e^{i\tau T_n}\big)&=\prod_{1\leqslant k\leqslant n}\Big(1+{\e^{i\tau k}-1\over k}\Big)\cr &=\exp\bigg\{S_n(\tau)+U(\tau)+W_n(\tau) +O\bigg({\tau\over  {n(1+n|\tau|)}}\bigg)\bigg\},\cr}\eqdef{carTn}$$
with  $$\eqalign{S_n(\tau)&:=\sum_{1\leqslant k\leqslant n}{\e^{i\tau k}-1\over k},\quad 
W_n(\tau):=\sum_{ k>n}{(\e^{i\tau k}-1)^2\over 2k^{2}}
\cr U(\tau)&:=\sum_{k\geqslant 1}\bigg\{\log \bigg(1+{\e^{i\tau k}-1\over k}\bigg)-{\e^{i\tau k}-1\over k}\bigg\}.\cr}$$
For  $|\tau|<2\pi$,  we may write
$$S_n(\tau)=\sum_{1\leqslant k\leqslant n}\int_0^{i\tau}\e^{kv}\d v=\int_0^{i\tau}{\e^{nv}-1\over 1-\e^{-v}}\d v=\int_0^n{\e^{i\tau v}-1\over v}\d v+V_n(\tau)$$
with $$\eqalign{V_n(\tau)&:=\int_0^{i\tau}(\e^{nv}-1)\Big({1\over 1-\e^{-v}}-{1\over v}\Big)\d v=\int_0^1(\e^{in\tau v}-1)g_\tau(v)\d v,\cr g_\tau(v)&:={i\tau\over 1-\e^{-i\tau v}}-{1\over v}\quad(0\leqslant v\leqslant 1).\cr}$$
Since $g_\tau(v)$ is twice continuously differentiable on $[0,1]$, partial integration yields
$$V_n(\tau)=V(\tau)+ {{a(\tau)\e^{in\tau}-\fs12\over n}}+O\Big({1\over n^2}\Big)\qquad (|\tau|\leqslant \pi),$$
with $$V(\tau):=-\int_0^1g_\tau(v)\d v,\quad a(\tau):={g_\tau(1)\over i\tau}={1\over 1-\e^{-i\tau}}-{1\over i\tau},$$ 
We have $W_n(\tau)\ll \tau$ if $|\tau|\leqslant 1/n$.
When $1/n\leqslant |\tau|\leqslant \pi $, we have
$$\eqalign{W_n(\tau)&=\sum_{ k> n}{ \e^{2ik\tau}-2\e^{i k\tau}\over 2k^2}+{1\over 2n} +O\Big({1\over n^2}\Big) ={1\over 2n}\bigg\{ 1+ O\Big({1\over 1+n\min(|\tau|,\pi-|\tau|)}
\Big)\bigg\}. \cr}\eqdef{majWn}$$
by Abel's summation. This estimate is hence also valid for $|\tau|\leqslant 1/n$, and so we deduce that
$$ V_n(\tau)+ W_n(\tau)=V(\tau) +{{a(\tau)\e^{in\tau} \over n}} +O\bigg({\tau\over  {n(1+n|\tau|)}}+{1\over n+n^2\min(|\tau|,\pi-|\tau|)}\bigg).$$
\goodbreak
 Put  $F(\tau):=\e^{U(\tau)+V(\tau)}-1$, so that $F(0)=0$ and $F$ may be analytically continued to the disc \hbox{$\{z\in\CC:|z|<2\pi\}$}.  We finally get 
$$
\EE(\e^{i\tau T_n})= {\wh{\varrho_0}(n\tau)}\bigg\{1+F(\tau)+W^*_n(\tau)+O\bigg({\tau\over 1+n^2\tau^2}\bigg)\bigg\}\qquad (n\geqslant 1,\, |\tau|\leqslant \pi ),\eqdef{EeitauTn}$$
with
$$W_n^*(\tau):=  {{a(\tau)\{1+F(\tau)\}\e^{in\tau} \over n}+O\bigg({1\over n+n^2\min(|\tau|,\pi-|\tau|)}\bigg) }\cdot\eqdef{EeitauTn+}$$
 
 \par 
It follows that, for $m\geqslant 1$, 
$$\eqalign{\PP&(T_n=m) ={1\over 2\pi}\int_{-\pi}^\pi\EE(\e^{i\tau T_n})\e^{-im\tau}\d\tau\cr
&={ 1\over 2\pi n}\int_{-n\pi}^{n\pi}\wh{\varrho_0}(\tau)\e^{-i\tau m/n}   \bigg\{1+F\Big({\tau\over n}\Big) +W_n^*\Big({\tau\over n}\Big)+O\bigg({\tau\over n(1+ \tau^2)}\bigg)\bigg\} \d\tau.\cr}\eqdef{ptnm}$$
\par 
By \eqref{tfr} and, say, \citeplus{Te15}{lemma III.5.9}, we have $$\wh{\varrho_0}(\tau)={-1\over i\tau}+O\Big({1\over \tau(1+|\tau|)}\Big)\quad (\tau\neq0),\qquad \wh{\varrho_0}(\tau)\asymp{1\over 1+|\tau|}\quad(\tau\in\r).\eqdef{whr0}$$
Therefore,   
the error term  of \eqref{ptnm}  contributes $\ll  1/n^2$ to the right-hand side. Summing over $m\leqslant  M(n) $, we obtain that the corresponding contribution to the right-hand side of \eqref{est1vn} is $\ll(\log n) /(n\log_2n)$, in accordance with  \eqref{encvn}.\par 
\par 
\par  We first evaluate 
$${1\over 2\pi n}\int_{-n\pi}^{n\pi}\wh{\varrho_0}(\tau)\e^{-i\tau m/n}\d\tau \qquad \big(n\geqslant 1,\,1\leqslant m\leqslant M(n)\big)\eqdef{mt}$$
by extending the integration range to $\r$ and inserting the first estimate~\eqref{whr0} to bound the integral over $\r\sset[-\pi n,\pi n]$. 
This yields
$$\eqalign{{1\over 2\pi n}\int_{-n\pi}^{n\pi}\wh{\varrho_0}(\tau)\e^{-i\tau m/n}\d\tau -{\varrho_0(m/n)\over n}
&=  { -1\over  \pi n}\int_{n\pi}^{\infty} {\sin(\tau m/n)\over \tau}\d\tau+O\Big({1\over n^2}\Big)
\cr&={(-1)^{m+1}\over   \pi^2  m n} +O\Big({1\over m^2n }+{1\over n^2}\Big).\cr}\eqdef{estmainrho}$$
\par 
In order to estimate  the contributions from $F$ and $W^*_n$  to the main term of \eqref{ptnm}, we use the more precise formula
$$\wh{\rho_0}(\tau)= {i\over \tau}-{\e^{i\tau}\over \tau^2}+O\Big({1\over \tau^3}\Big) \quad(|\tau|\geqslant 1).\eqdef{whr0+}$$
\par Writing $F(\tau)=\tau G(\tau)$, we indeed  deduce from by \eqref{whr0+}  that 
$$\eqalign{\int_{-n\pi}^{n\pi}\wh{\varrho_0} (\tau)\e^{-i\tau m/n} F\Big({\tau\over n}\Big) &\d\tau=\int_{-n\pi}^{n\pi}\tau \wh{\varrho_0}(\tau){\e^{-i\tau m/n}\over n} G\Big({\tau\over n}\Big)\d \tau \cr&=i\int_{I_n}  {\e^{-i\tau m } } G( \tau )\Big(1+i{\e^{i\tau n}\over n\tau }\Big)\d \tau+O\Big({1\over n}\Big),\cr}$$
with $I_n:=[-\pi,\pi]\sset[ -1/n, 1/n ]$. A standard computation furnishes  $G(\pi)-G(-\pi)=-2/\pi$.\note{ $V(\pm\pi)=\log (\pi/2)\mp\fs12 i\pi$, $\e^{U(\pm\pi)}=-2\e^{\gamma}$, $F(\pm\pi)= \mp i\pi \e^\gamma-1$. 
} Integrating by parts, we get 
$$\eqalign{i\int_{I_n}  {\e^{-i\tau m } } G( \tau ) \d \tau
\ &=2 {(-1)^{m+1}\over\pi m} +{1\over m}\int_{- \pi}^{ \pi} {\e^{-i\tau m } } G'( \tau )\d \tau+O\Big( {1\over n}\Big)
\cr&=2 {(-1)^{m+1}\over \pi m}+O\Big({1\over m^2}+{1\over n}\Big),\cr}$$
and, similarly, $${-1\over n}\int_{I_n}  {\e^{-i\tau m } } G( \tau ) {\e^{i\tau n}\over \tau } \d \tau={-1\over n } \int_{I_n}     {G( \tau ) -G(0)\over  \tau }{\e^{i\tau (n-m)}  } \d \tau +O\Big({1\over n}\Big) \ll {1\over n}\cdot $$
  We can thus state that, for $n\geqslant 1,\,1\leqslant m\leqslant M(n)$, we have
$$ {1\over 2\pi n}\int_{-n\pi}^{n\pi}\wh{\varrho_0} (\tau)\e^{-i\tau m/n} F\Big({\tau\over n}\Big) \d\tau={(-1)^{m +1 }\over  \pi^2mn}+O\Big({1\over m^2n}+{1\over n^2}\Big).\eqdef{contF}$$ 
\par 
It remains to estimate the contribution to \eqref{ptnm} involving  $W_n^*(\tau/n) $. 
 The error term of \eqref{EeitauTn+} clearly contributes   $\ll 1/n^2$. Arguing as for the proof of \eqref{estmainrho}, we finally show that the contribution to \eqref{ptnm} arising from $  {a(\tau/n)(1+F(\tau/n))\e^{i \tau}/n} $ is also $\ll 1/n^2$.  

 The above estimates and \eqref{estmainrho} furnish together 
$$ \PP(T_n=m)-{1\over n}\varrho_0\Big({m\over n}\Big)=(-1)^{m+1} { 2 \over \pi^2 m n} +O\Big({1\over m^2n }+{1\over n^2}\Big)\qquad \big(1\leqslant m\leqslant M(n)\big).\eqdef{cal}$$
Summing over $m\leqslant M(n)$ provides the required estimate \eqref{encvn}.

\paraun{Proof of \ref{T2}}
We first note that \eqref{cal} Êyields 
$$\PP(T_n=m)={1\over n}\varrho_0\Big({m\over n}\Big)+O\Big({1\over n^2}\Big)\qquad (m\asymp n\geqslant 1).\eqdef{PTnm}$$
Hence, writing $L_N$ for $L_N(u)$ here and throughout,
$$\EE(L_N)=\sum_{n\leqslant N}{1\over n}\varrho_0\Big({m_n\over n}\Big)+O(1)=\varrho_0(u)\log N+O(\eta_N\log N+1)\qquad (N\geqslant 1).\eqdef{ELN}$$
The stated result will follow from an estimate of the variance $\VV(L_N)$. We have
$$\EE(L_N^2)=\EE(L_N)+2\sum_{1\leqslant \nu<n\leqslant N}\PP(T_\nu=m_\nu)\alpha_{\nu n},\eqdef{ELN2}$$
with $$\alpha_{\nu n}:=\PP\big(T_n-T_\nu=m_n-m_\nu\big)\qquad (1\leqslant \nu<n\leqslant N).$$
\goodbreak
Let us initially assume that $\{m_n\}_{n\geqslant 1}$ is strictly increasing and hence that $u\geqslant 1$. We note right away that, for large $n$ and  $\nu>(u+1)n/(u+2)$, we have $\alpha_{\nu n}=0$,  since the corresponding event is then impossible: either $T_n-T_\nu>\nu>(u+1)(n-\nu)>m_n-m_\nu$ or $T_n-T_\nu=0\neq m_n-m_\nu$. Therefore, we may assume $\nu\leqslant n(u+1)/(u+2)$ in the sequel.\par 
\smallskip

By \eqref{carTn}, we have, writing $\varphi_j(\tau):=\EE\big(\e^{i\tau T_j}\big)$ $(j\geqslant 1)$,
$$\alpha_{\nu n}={1\over 2\pi}\int_{-\pi}^\pi{\varphi_n(\tau)\over \varphi_\nu(\tau)}\e^{-i(m_n-m_\nu)\tau}\d \tau=\beta_{\nu n}+\Delta_{\nu n},\eqdef{abd}$$
where $\beta_{\nu n}$ denotes the contribution of the interval $|\tau|\leqslant 1/2\nu$, and $\Delta_{\nu n}$ that of the complementary range $1/2\nu<|\tau|\leqslant \pi$.
\par 
Now we may derive from \eqref{EeitauTn} that
$$\beta_{\nu n}={1\over 2\pi n}\int_{-n/2\nu}^{n/2\nu}{\wh{\varrho_0}(\tau)\over \varphi_\nu(\tau/n)}\Big\{1+F\Big({\tau\over n}\Big)+O\Big({1\over n}\Big)\Big\}\e^{-i(m_n-m_\nu)\tau/n}\d \tau.$$
Invoking  \eqref{whr0} and \eqref{EeitauTn}-\eqref{EeitauTn+} with $\nu$ in place of $n$ to estimate the contribution of the error term, we get
$$\beta_{\nu n}=\beta^*_{\nu n}+O\bigg({\log n\over n^2}\bigg)\eqdef{aa*}$$
with 
$$\beta^*_{\nu n}:={1\over 2\pi n}\int_{-n/2\nu}^{n/2\nu}{\wh{\varrho_0}(\tau)\over \varphi_\nu(\tau/n)}\Big(1+F\Big({\tau\over n}\Big)\Big)
\e^{-i(m_n-m_\nu)\tau/n}\d \tau.\eqdef{a*1}$$
Note that the remainder term of \eqref{aa*} in turn contributes $\ll  1$ to \eqref{ELN2}. 
\par  
Still assuming that $\nu\leqslant (u+1)n/(u+2)$ and observing that $$\psi_\nu(z):=\prod_{k\leqslant \nu}\{1+(\e^{kz}-1)/k\}^{-1}\ll 1\qquad (z\in\CC,\,|z|\leqslant 2/3\nu),$$  we may write
$${1+F(\tau/n)\over \varphi_\nu(\tau/n)}=1+\sum_{j\geqslant 1}\mu_{\nu j}\Big({\tau\over n}\Big)^{j}\qquad (|\tau|\leqslant n/2\nu)$$
with $\mu_{\nu j}\ll (3\nu/2)^j$ $(j\geqslant 1)$. 
Hence,  in view of \eqref{whr0}, the contribution of the series to \eqref{a*1}~is 
$$\eqalign{&{1\over 2\pi n}\sum_{j\geqslant 1}{\mu_{\nu j}\over n^j}\int_{-n/2\nu}^{ n/2\nu}\bigg\{\tau^{j-1}\e^{i\tau (m_n-m_\nu)/n}+O\bigg({\tau^{j-1}\over 1+|\tau| }\bigg)\bigg\}\d \tau\cr&
\ll{1\over n}\sum_{j\geqslant 1}\mu_{\nu j}\int_{-1/ 2\nu}^{1/ 2\nu }\bigg\{v^{j-1}\e^{iv(m_n-m_\nu)}+O\bigg({v^{j-1}\over 1+n|v|} \bigg)\bigg\}\d v
\cr&\ll {\nu\over n(n-\nu)}+{\nu\over n^2} +{\nu\log (2n/\nu)\over n^2}\ll {\nu\log (2n/\nu)\over n^2}
,\cr}$$
since the assumption that $\{m_n\}_{n=1}^{\infty}$ is strictly increasing implies $m_n-m_\nu\geqslant n-\nu$.
\par 
Arguing as in the proof of \eqref{mt} to evaluate the main term, we eventually get,  
$$\beta^*_{\nu n}={1\over n}\varrho_0\Big({m_n-m_\nu\over n}\Big)+O\Big({ {\nu\log (2n/\nu)}\over n^2}\Big)\qquad (1\leqslant \nu<n).\eqdef{estbeta*}$$
\par \goodbreak
We now consider $\Delta_{\nu n}$.  Let $$S_{\nu n}(\tau):=\sum_{\nu<k\leqslant n}{\e^{ik\tau}\over k},$$
and note right away that $S_{\nu n}(\tau)$ is bounded for $1/2\nu<|\tau|\leqslant \pi$.
For $\nu\geqslant 1$, we have
$${\varphi_n(\tau)\over \varphi_\nu(\tau)}={\nu\over n}\prod_{\nu<k\leqslant n}\Big(1+{\e^{ik\tau}\over k-1}\Big)={\nu\over n}\exp\bigg\{S_{\nu n}(\tau)+O\Big({1\over \nu+|\tau|\nu^2}+{1\over \nu+(\pi- {|\tau|})\nu^2}\Big)\bigg\},$$
where we used the bounds
$$\eqalign{&\sum_{\nu<k\leqslant n}{\e^{ ik\tau}\over k(k-1)}\ll{1\over \nu+|\tau|\nu^2},
\quad\sum_{\nu<k\leqslant n}{\e^{2ik\tau}\over (k-1)^2}\ll{1\over \nu+|\tau|\nu^2}+{1\over \nu+(\pi- |\tau| )\nu^2},\cr
&\sum_{\nu<k\leqslant n}{\e^{ikj\tau}\over (k-1)^j}\ll{1\over \nu^{j-1}}\quad(j> 2) .\cr}$$
Carrying back into \eqref{abd}, we obtain
$$\Delta_{\nu n}={\nu\over \pi n}\Re \int_{1/2\nu}^{\pi}\e^{S_{\nu n}(\tau)-i(m_n-m_\nu)\tau}\d\tau+O\Big({\log 2\nu\over n\nu}\Big).\eqdef{Dnnu}$$ 
Now observe that $|S_{\nu n}'(\tau)|\leqslant \pi/|\tau|\leqslant \dm( n-\nu)\leqslant  \dm(m_n-m_\nu)$ if, say $\nu\leqslant n/15$ or $|\tau|>10(u+1)/\nu$. Hence, on this assumption, a standard estimate on trigonometric integrals such as \citeplus{Ti51}{Lemma 4.2} furnishes the bound $\ll1/(m_n-m_\nu)\ll1/(n-\nu)\ll1/n$ for the last integral. Since, in the case $\nu>n/15$, the contribution of the range $1/2\nu<|\tau|\leqslant 10(u+1)/\nu$ to the same integral is trivially $\ll 1/\nu$, we finally get 
$$\Delta_{\nu n}\ll {\log 2\nu\over n\nu}+{\nu\over n^2}\qquad \big(1\leqslant \nu\leqslant (u+1)n/(u+2)\big).\eqdef{estDnun}$$
 \par 
Gathering our estimates and using the fact that $\varrho$ is Lipschitz on $[0,\infty[$, we obtain
$$\alpha_{\nu n}={1\over n}\varrho_0\Big({m_n\over n}\Big)+O\Big({\nu\log (2n/\nu) \over n^2}+{\log 2\nu\over n\nu}\Big)\quad(1\leqslant \nu<n).\eqdef{evalphanun}$$
Carrying back into \eqref{ELN2} and applying \eqref{PTnm} for the pair $(\nu,m_\nu)$ yields
$$\EE(L_N^2)-\EE(L_N)=2\sum_{1\leqslant \nu<n\leqslant N}\bigg\{\varrho_0\Big({m_n\over n}\Big)\varrho_0\Big({m_\nu\over \nu}\Big){1\over \nu n}+O\bigg({ \log (2n/\nu)\over n^2}+{\log 2\nu\over n\nu^2}\bigg)\bigg\},$$
and hence
$$\VV(L_N)\ll\log N.$$
 Selecting $N=N_k:=2^{k^3}$  for $k\geqslant 1$, we deduce from the Borel-Cantelli lemma that, given any $\varepsilon>0$, the estimate
$$L_{N_k}-\EE(L_{N_k})\ll k^2(\log 2k)^{1/2+\varepsilon} $$
holds almost surely. In view of \eqref{ELN}, this implies the stated result since $L_N$ is a non-decreasing function of $N$.
\par 
We next consider the case of a non-decreasing sequence $\{m_n\}_{n\geqslant 1}$. Accordingly, we fix $u>0$. By hypothesis, for some integer $q=q_N\geqslant 2$ such that $q_N=o\big(\log N\big)$ as $N\to\infty$, we have $m_n>m_\nu$ whenever $n-\nu\geqslant q_N$.\par 
 Put
$$L_N(u;a):=\sum_{\tri{n\leqslant N}{n\md aq}{T_n=m_n}}1\quad(1\leqslant a\leqslant q).$$
By \eqref{PTnm}, we have, for all $a\in[1,q]$,
$$\EE\big(L_N(u;a)\big)=\sum_{\di{n\leqslant N}{n\md aq}}{1\over n}\varrho_0\Big({m_n\over n}\Big)+O\Big({1\over a^2}\Big),$$
and, by \eqref{evalphanun},
$$\eqalign{\VV\big(L_N(u;a)\big)-\EE\big(L_N(u;a)\big)&\ll\sum_{\di{1\leqslant \nu<n\leqslant N}{\nu,n\md aq}}\Big\{{\log 2n/\nu\over n^2}+{\log 2\nu\over n\nu^2}\Big\}
\cr&\ll\sum_{\di{q<n\leqslant N}{n\md aq}}\Big\{{1\over nq}+{\log 2a\over a^2n}\Big\}\ll{\log N\over q}\Big\{{1\over q}+{\log 2a\over a^2}\Big\}.\cr}$$ 
Summing over $a\in[1,q]$, we get
$$\VV\big(L_N(u)\big)\leqslant q\sum_{1\leqslant a\leqslant q}\VV\big(L_N(u;a)\big)\ll q\EE(L_N(u))+\log N\ll q\log N.$$
This is all needed.
\par 
\medskip
 \noi {\it Acknowledgment.} The authors express warm thanks to Eugenijus Manstavi\v cius for precise and useful comments on a preliminary version of the manuscript. 
 \goodbreak  
\bigskip\bigskip
\centerline{\twelvebf References}
\bigskip
\eightpoint{
 \bibtem{ABT93} R. Arratia, A.D. Barbour, S. TavarŽ, On random polynomials over finite fields, {\it Math. Proc. Cambridge Philos. Soc. \bf114} (1993), no. 2, 347Ð368. \par 
 \bibtem{ABT99} R. Arratia, A.D. Barbour, S. TavarŽ, The Poisson-Dirichlet distribution and the scale-invariant Poisson process, {\it Combin. Probab. Comput. \bf8} (1999), no. 5, 407Ð416.\par 
 \bibtem{ABT00} R. Arratia, A.D. Barbour, S. TavarŽ, Limits of logarithmic combinatorial structures, {\it Ann. Probab. \bf28} (2000), no. 4, 1620Ð1644. 
 \par 
 \bibtem{BG19} C. Bhattacharjee \& L. Goldstein, Dickman approximation in simulation, summations and perpetuities, {\it Bernoulli \bf 25} (2019), no. 4A, 2758Ð2792. \par 
 \bibtem{Bi72}  P. Billingsley,  On the distribution of large prime factors, {\it Period. Math. Hungar. \bf2}
(1972), 283--289.\par 
\bibtem{Ca87} M. Car, ThŽormes de densitŽ dans $\F_q$[X], {\it Acta Arith. \bf48} (1987), no. 2, 145Ð165.
\par
\bibtem{CH03}  W.-M. Chen \& H.-K. Hwang, Analysis in distribution of two randomized algorithms for finding the maximum in a broadcast communication model, {\it J. Algorithms \bf46} (2003), no. 2, 140Ð177.
\par 
\bibtem{De01} L. Devroye, Simulating perpetuities, {\it Meth. Comput. Appl. Probab. \bf3}  (2001), 97--115.
\par
 \bibtem{DF10} L. Devroye \& O. Fawzi, Simulating the Dickman distribution, {\it Statist. Probab. Lett. \bf80} (2010), no. 3-4, 242Ð247.
 \par 
\bibtem{Di30} K. Dickman, On the frequency of numbers containing primes of a
 certain relative magnitude, {\it Ark. Mat. Astr. Fys. \bf 22} (1930), 1--14.\par 
 \bibtem{El79-80} P.D.T.A. Elliott, Probabilistic number theory : mean value theorems, Grund\-lehren der Math. Wiss.  239; Probabilistic number theory : central limit theorems, {\it ibid.} 240, Springer-Verlag, New York, Berlin, Heidelberg 1979, 1980. 
 \bibtem{FH10}  J.A. Fill \& M.L. Huber, Perfect simulation of Vervaat perpetuities {\it Electron. J. Probab. \bf15} (2010), no. 4, 96Ð109. 
 \par   
\bibtem{GSW18} R. Giuliano, Z.S. Szewczak \& M.J.G. Weber, Almost sure local limit theorem for the Dickman
distribution, {\it Period. Math. Hungar.  \bf76}, no. 2 (2018),  155Ð197. 
\par 
\bibtem{GG96}  C.M. Goldie \& R. GrŸbel, Perpetuities with thin tails, {\it Adv. in Appl. Probab. \bf28} (1996), no. 2, 463Ð480.
\par 
\bibtem{HT02} H.K. Hwang \& T.-H. Tsai, Quickselect and the Dickman Function, {\it Combinatorics, Probability and Computing \bf11}, (2002) 353--371.
\par 
\bibtem{Ki77} J.F.C. Kingman, The population structure associated with the Ewens sampling formula,
{\it Theoretical Population Biology \bf11} (1977), 274Ð283.
\par 
\bibtem{KM97}  A. Knopfmakher \& E. Manstavi\v cius, On the largest degree of an irreducible factor of a polynomial in $\F_q[X]$ (Russian), {\it Liet. Mat. Rink. \bf37} (1997), no. 1, 50Ð60; translation in: {\it Lithuanian Math. J. \bf37} (1997), no. 1, 38Ð45.
\par 
\bibtem{Ma92} E. Manstavi\v cius, Remarks on elements of semigroups that are free of large prime factors (Russian), {\it Liet. Mat. Rink. \bf32} (1992), no. 4, 512Ð525; translation in: {\it Lithuanian Math. J. \bf32} (1992), no. 4, 400Ð409 (1993).\par 
\bibtem{MP17} E. Manstavi\v cius \& R. Petuchovas, Local probabilities and total variation distance for random permutations, {\it Ramanujan J. \bf43} (2017), no. 3, 679Ð696.
\par 
\bibtem{Pi19}  R. G. Pinsky, A natural probabilistic model on the integers and its relation to Dickman-type distributions and Buchstab's function, in: {\it Probability and analysis in interacting physical systems}, 267Ð294, Springer Proc. Math. Stat. 283, Springer, Cham, 2019.
\par 
\bibtem{SL66}  L.A. Shepp \& S.P. Lloyd, Ordered cycle lengths in a random permutation, {\it Trans. Amer. Math. Soc. \bf121} (1966), 340--357.
\par 
\bibtem{Te99} G. Tenenbaum, Crible d'\'Eratosthne et modle de Kubilius, in: K. Gy\H ory, H. Iwaniec, J. Urbanowicz
(eds.), {\it Number Theory in Progress}, Proceedings of the conference in honor of Andrzej Schinzel,
Zakopane, Poland 1997, 1099--1129, Walter de Gruyter, Berlin, New York.
\par 
\bibtem{Te00} G. Tenenbaum, A rate estimate in Billingsley's theorem for the size distribution of large prime
factors, {\it Quart. J. Math. (Oxford) \bf51} (2000), 385--403.
\bibtem{Te15} G. Tenenbaum, {\it Introduction to analytic and probabilistic number theory}, 3rd ed., Graduate Studies in Mathematics 163, Amer. Math. Soc. 2015.\par
\bibtem{Ti51} E.C. Titchmarsh,
 {\it The theory of the Riemann zeta-function}, Oxford University Press, 1951.
 \par 
}

{\leftskip2mm\rightskip-2cm\sevenrm
\gutter=4cm \doublecolumns
 \obeylines \baselineskip=7pt
RŽgis de la Bretche
UniversitŽ de Paris, Sorbonne UniversitŽ, CNRS, 
Institut de Math. de Jussieu-Paris Rive Gauche,
 F-75013 Paris,  
France
\smallskip
{\seventt regis.de-la-breteche@imj-prg.fr}
 GŽrald Tenenbaum\par
Institut \'Elie Cartan\par 
Universit\'e de Lorraine\par
 BP 70239\par
54506 Vand\oe uvre-ls-Nancy Cedex\par
 France
\smallskip
{\seventt gerald.tenenbaum@univ-lorraine.fr}
\singlecolumn
\par}
\end

%% file: itenn.tex


  %
  %

  %
  %

  \catcode`@=12 

 \def\defrefnote#1{\definexref{#1}{{\the\footnotenumber}}{refnotes}}

  %
  %


\ifx\couleurs\oui
\input graphicx
\fi
\input eplain.tex
 
\ifx\couleurs\oui
\beginpackages
\usepackage{color}
\endpackages 
\long\def\rge#1{{\color{red}#1}}

\definecolor{bleu-iecn}{cmyk}{.98,.13,.1,.55}

\fi

\makeatletter
\def\numberedfootnote{%
ÊÊ\global\advance\footnotenumber by 1
ÊÊ\@eplainfootnote{{\number\footnotenumber}}%
}%
\def\makecolumns#1/#2 {\par \begingroup
ÊÊ \@columndepth = #1
ÊÊ \advance\@columndepth by -1
ÊÊ \divide \@columndepth by #2
ÊÊ \advance\@columndepth by 1
ÊÊ \@linestogoincolumn = \@columndepth
ÊÊ \@linestogo = #1
ÊÊ \currentcolumn = 1
ÊÊ \def\@endcolumnactions{%
ÊÊÊÊÊÊ\ifnum \@linestogo<2
ÊÊÊÊÊÊÊÊ \the\crtok \egroup \endgroup \par 
ÊÊÊÊÊÊ\else
ÊÊÊÊÊÊÊÊ \global\advance\@linestogo by -1
ÊÊÊÊÊÊÊÊ \ifnum\@linestogoincolumn<2
ÊÊÊÊÊÊÊÊÊÊÊÊ\global\advance\currentcolumn by 1
ÊÊÊÊÊÊÊÊÊÊÊÊ\global\@linestogoincolumn = \@columndepth
ÊÊÊÊÊÊÊÊÊÊÊÊ\the\crtok
ÊÊÊÊÊÊÊÊ \else
ÊÊÊÊÊÊÊÊÊÊÊÊ&\global\advance\@linestogoincolumn by -1
ÊÊÊÊÊÊÊÊ \fi
ÊÊÊÊÊÊ\fi
ÊÊ }%
ÊÊ \makeactive\^^M
ÊÊ \letreturn \@endcolumnactions
ÊÊ \@columnwidth = \hsize
ÊÊÊÊ \advance\@columnwidth by -\parindent
ÊÊÊÊ \divide\@columnwidth by #2
ÊÊ \penalty\abovecolumnspenalty
ÊÊ \noindent 
ÊÊ \valign\bgroup
ÊÊÊÊ &\hbox to \@columnwidth{\strut \hsize = \@columnwidth ##\hfil}\cr
}%
\makeatother

\lefteqnumbers
   \def\testd{oui}
   \def\choixlat{\ifx\numadroite\testd\righteqnumbers
            \else  \lefteqnumbers\fi}
    \choixlat

\catcode`@=\letter
\def\@eplainfootnote#1{\let\@sf\empty 
  \ifhmode\edef\@sf{\spacefactor\the\spacefactor}\/\fi
  \global\advance\hlfootlabelnumber by 1
  \hlstart@impl{foot}{\hlfootlabel}%
  \hldest@impl{footback}{\hlfootbacklabel}%
  \hbox{$^{(#1)}$}%
  \hlend@impl{foot}%
  \@sf\vfootnote{#1.}%
}%
\catcode`@=\other

  \interfootnoteskip=0pt
  \let\note=\numberedfootnote
  \everyfootnote={\eightpoint\leftskip=5truemm\rightskip5truemm}
  
  \hsize150truemm\vsize 240truemm\hoffset=5truemm
  \def\dimstand{\hsize 150truemm\vsize 240truemm\hoffset=5truemm\voffset=0truemm}
  
  \pretolerance=500\tolerance=1000\brokenpenalty=5000
  \parindent3mm
  
  \countdef\temps=170
  \temps=\time
  \countdef\nminutes=171{\nminutes = \time}
  \countdef\nheure=172
  \def\heure{\begingroup                     
     \temps = \time \divide\temps by 60
     \nheure = \temps                        
     \nminutes = \time
     \multiply\temps by 60
     \advance\nminutes by -\temps            
     \ifnum\nminutes<10 \toks1 = {0}%
     \else\toks1 = {}%
     \fi
     \number\nheure h\the\toks1 \number\nminutes  
  \endgroup}%

  \newcount\chstart
  \chstart=\pageno
 \headline={\ifnum\pageno=\chstart {\hfill} \else {\hss \tenrm --\ \folio\ --\hss}\fi}
  \footline={\hfill}
  \normalbaselines
  \frenchspacing
    \def\dater{\vglue-10mm\rightline{(\the\day/\the\month/\the\year)}}
  \def\dateheure{\vglue-10mm\rightline{(\the\day/\the\month/\the\year,\ \heure)}}

  \newif\ifpagetitre \pagetitretrue
\newtoks\hautpagetitre \hautpagetitre={\hfill}
\newtoks\baspagetitre \baspagetitre={\hfill}
\newtoks\auteurcourant \auteurcourant={\hfill}
\newtoks\titrecourant \titrecourant={\hfill}
\newtoks\hautpagegauche
\newtoks\hautpagedroite
\newtoks\hautpagemilieu
\hautpagemilieu={\tenrm\hfil -- \folio\ -- \hfil}
\hautpagegauche={\ifx\midfolio\oui\the\hautpagemilieu\else\tenrm\folio\hfill\the\auteurcourant\hfill\fi}
\hautpagedroite={\ifx\midfolio\oui\the\hautpagemilieu\else\hfill\the\titrecourant\hfill\tenrm\folio\fi}
\newtoks\baspagegauche \baspagegauche={\hfil}
\newtoks\baspagedroite \baspagedroite={\hfil}
\headline={\ifpagetitre\the\hautpagetitre
\else\ifodd\pageno\the\hautpagedroite\else\the\hautpagegauche\fi\fi }
\footline={\ifpagetitre\the\baspagetitre
\else\ifodd\pageno\the\baspagedroite
\else\the\baspagegauche\fi\fi \global\pagetitrefalse}

\def\pageblanche{\vfill\eject\pagetitretrue
\null\vfill\eject
\pagetitretrue
}
\def\chgtpage{\ifodd\pageno \else
\pageblanche \fi \pagetitretrue\titreun=0\footnotenumber=0}

\def\chgtpageincrtitreun{\ifodd\pageno \else
\pageblanche \fi \pagetitretrue\footnotenumber=0}

\def\majnombres{\ifodd\pageno \else
\pageblanche \fi \pagetitretrue\hautpoly\titreun=0\footnotenumber=0}

\def\hautspages#1#2{\auteurcourant={\ninepcap#1}\titrecourant={\nineit#2}}

\ifnum\chstart=\pageno \pagetitretrue\fi
  

  \def\PAR{\par}
  
  \def\leftnote#1{\vadjust{\setbox1=\vtop{\hsize 20mm\parindent=0pt\eightpoint
  \baselineskip=9pt\rightskip=4mm plus 4mm\vskip-4.7mm#1}\hbox{\kern-2cm\smash{\box1}}}}

  
  \def\raggedcenter{\leftskip=20pt plus 10em  
       \rightskip=\leftskip 
        \parfillskip=0pt 
         \spaceskip=.3333em \xspaceskip=.5em 
          \pretolerance=9999 \tolerance=9999
           \hyphenpenalty=9999 \exhyphenpenalty=9999 }
           
  \def\titrecentre#1{{\parindent0mm\raggedcenter
       \spaceskip=.6em plus .2em minus .2em\xspaceskip=.6em plus .2em minus .2em
        \tit#1\par}}
        


  \def\oui{oui}
  
   \def\fontetitreun{\ifx\paradouze\oui\douzepts\gpdouze\twelvebf\textfont1=\twelveib\else
\quatorzepts\gpquatorze\fourteenbf\fi}

\def\fontetitreunl{\douzepts\textfont1=\twelveib\scriptfont1=\tenib\fourteenti}
 
 \def\fontetitredeux{\textfont1=\eleveni\ifx\paradouze\oui\onzepts\scriptfont1=\ninei\elevenit\else
                        \douzepts\twelveit\fi}
 
   \def\fontetitredeuxb{\ifx\paradouze\oui\onzepts\eleventi\gponze\textfont1=\elevenib\scriptfont1=\nineib
                         \else\douzepts\twelveti\scriptfont1=\twelveib\scriptfont1=\tenib\gpdouze\fi}
                         
\def\fontetitredeuxl{\onzepts\textfont1=\elevenbf\scriptfont1=\ninebf\twelvebf}
  
\def\fontetitretrois{\textfont0=\elevenrm\scriptfont0=\eightrm\textfont1=\eleveni
                      \scriptfont1=\eighti\scriptscriptfont1=\sixi\elevenit}
                      
\def\fontetitrequatre{\textfont0=\elevenrm\scriptfont0=\eightrm\textfont1=\eleveni
                      \scriptfont1=\eighti\scriptscriptfont1=\sixi\elevenrm}
  
  \newcount\titreun\titreun=0
  \newcount\titredeux\titredeux=0
  \newcount\titretrois\titretrois=0
  \newcount\titrequatre\titrequatre=0
  \newcount\enonce\enonce=0
  
  \def\incr#1{\global\advance#1 by 1 {\the #1}}
  \def\avance#1{\global\advance#1 by 1}
  \def\init#1{\global#1=0}
  
  \long\def\Indentation#1#2{\setbox10=\hbox{\fontetitreun#1}
  	                    \ifdim\wd10 < 4mm
                         \setbox10=\hbox to 4mm{\box10\hfill}
                       \else\ifdim\wd10 < 6mm
                         \setbox10=\hbox to 6mm{\box10\hfill}
  	                    \else\ifdim\wd10 < 8mm
                         \setbox10=\hbox to 8mm{\box10\hfill}
                       \else\ifdim\wd10 < 12mm
                         \setbox10=\hbox to 12mm{\box10\hfill}
                       \fi\fi\fi\fi
                       \dimen10=\hsize
                       \advance \dimen10 by -\wd10
                       \noindent \box10 %
                       \ignorespaces
                       \hbox{\vtop{\hsize=\dimen10\raggedright\noindent\fontetitreun#2}}}

  \long\def\paraun#1{\removelastskip\par\medskip\goodbreak\vskip0pt plus.01\vsize\penalty-100
                \vskip0pt plus-.01\vsize
  	              \init{\titredeux}\ifnum\optionparag=1{\init\eqnumber\init\enonce}\else{}\fi
                  \goodbreak{\fontetitreun
  	                \Indentation{\incr{\titreun}.\ }{\fontetitreun #1\par}}\nobreak\medskip}

 %
 %
 \long\def\paraunc#1{\removelastskip\par\bigskip\goodbreak\vskip0pt plus.01\vsize\penalty-100
                \vskip0pt plus-.01\vsize
  	              \init{\titredeux}
                 \ifnum\optionparag=1{\init{\eqnumber}\init\enonce}\else{}\fi
                  \goodbreak
  	                {\parindent0mm\raggedcenter\fontetitreun\incr{\titreun}.\ 
                     \fontetitreun #1\par}\nobreak\medskip}
                     
\newtoks\titreunl
\titreunl={\ifnum\titreun=1{I}\fi%
\ifnum\titreun=2{II}\fi%
\ifnum\titreun=3{III}\fi%
\ifnum\titreun=4{IV}\fi%
\ifnum\titreun=5{V}\fi%
\ifnum\titreun=6{VI}\fi%
\ifnum\titreun=7{VII}\fi%
\ifnum\titreun=8{VIII}\fi%
\ifnum\titreun=9{IX}\fi%
\ifnum\titreun=10{X}\fi%
\ifnum\titreun=11{XI}\fi%
\ifnum\titreun=12{XII}\fi%
\ifnum\titreun=13{XIII}\fi%
}
\long\def\paraunl#1{\removelastskip\par\bigskip\bigskip\goodbreak\vskip0pt plus.01\vsize\penalty-100
                \vskip0pt plus-.01\vsize
  	              \init{\titredeux}\ifnum\optionparag=1{\init\eqnumber\init\enonce}\else{}\fi
                  \goodbreak{\fontetitreunl
  	                \Indentation{\global\advance\titreun by 1{\the\titreunl}.\ }{\fontetitreunl #1\par}}\nobreak\smallskip}

  
  \long\def\paradeux#1{\init{\titretrois}\vskip0pt plus.01\vsize\penalty-10
                \vskip0pt plus-.01\vsize\ifx \elie\oui\medskip\ifnum\titredeux>0\medskip\fi\fi
                 \Indentation{\fontetitredeux\the\titreun${\cdot}$\incr{\titredeux}.}
                              {\fontetitredeux\textfont1=\eleveni#1}\nobreak\par }
  
  \long\def\paradeuxb#1{\init{\titretrois}\vskip0pt plus.001\vsize\penalty-10
                \vskip0pt plus-.01\vsize{\ifx \elie\oui\medskip\ifnum\titredeux>0\medskip\fi\fi
                  \Indentation
  {\fontetitredeuxb\the\titreun${\cdot}$\incr{\titredeux}.}{ \fontetitredeuxb#1}}\nobreak
\smallskip}

\newtoks\titredeuxl
\titredeuxl={\ifnum\titredeux=1{A}\fi%
\ifnum\titredeux=2{B}\fi%
\ifnum\titredeux=3{C}\fi%
\ifnum\titredeux=4{D}\fi%
\ifnum\titredeux=5{E}\fi%
\ifnum\titredeux=6{F}\fi%
\ifnum\titredeux=7{G}\fi%
\ifnum\titredeux=8{H}\fi%
\ifnum\titredeux=9{I}\fi%
\ifnum\titredeux=10{J}\fi%
\ifnum\titredeux=11{K}\fi%
\ifnum\titredeux=12{L}\fi%
\ifnum\titredeux=13{M}\fi%
}
 \long\def\paradeuxl#1{\init{\titretrois}\vskip0pt plus.001\vsize\penalty-10
                \vskip0pt plus-.01
                \vsize \bigskip%
                  \Indentation
     {\fontetitredeuxl\global\advance\titredeux by 1
  \quad \the\titreunl${\cdot}$\the\titredeuxl.}{ \fontetitredeuxl#1}
  \removelastskip\nobreak\smallskip}
  

  \long\def\paratrois#1{\init{\titrequatre}\ifdim\lastskip<\smallskipamount
                \removelastskip\smallskip\fi
                 \vskip0pt plus.01\vsize\penalty-10
                  \vskip0pt
plus-.01\vsize{\ifx \elie\oui\ifnum\titretrois>0\medskip\fi\fi
\Indentation{\fontetitretrois\the\titreun${\cdot}$\the\titredeux${\cdot}$\incr{\titretrois}.\ }
  {\hskip0mm\baselineskip=14pt\fontetitretrois#1}\nobreak\smallskip}}
  
  
  \long\def\paratroisl#1{\init{\titrequatre}\ifdim\lastskip<\smallskipamount
                \removelastskip\fi
                 \vskip0pt plus.01\vsize\penalty-10
                  \vskip0pt
plus-.01\vsize\ifx \elie\oui\bigskip
\fi
\Indentation{\fontetitretrois\quad \quad \the\titreunl{${\cdot}$}\the\titredeuxl${\cdot}$\incr{\titretrois}.\ }
  {\hskip0mm\fontetitretrois#1}\nobreak\smallskip}


  \long\def\paraquatre#1{\ifdim\lastskip<\smallskipamount
                \removelastskip\smallskip\fi
                 \vskip0pt plus.01\vsize\penalty-10
                  \vskip0pt
                  plus-.01\vsize\par
 
\Indentation{\fontetitrequatre \the\titreun{${\cdot}$}\the\titredeux${\cdot}$\the\titretrois${\cdot}$\incr{\titrequatre}.\ }
{\hskip0mm\fontetitrequatre#1}\nobreak\smallskip}


\newtoks\titrequatrel
\titrequatrel={\ifnum\titrequatre=1{a}\fi%
\ifnum\titrequatre=2{b}\fi%
\ifnum\titrequatre=3{c}\fi%
\ifnum\titrequatre=4{d}\fi%
\ifnum\titrequatre=5{e}\fi%
\ifnum\titrequatre=6{f}\fi%
\ifnum\titrequatre=7{g}\fi%
\ifnum\titrequatre=8{h}\fi%
\ifnum\titrequatre=9{i}\fi%
\ifnum\titrequatre=10{j}\fi%
\ifnum\titrequatre=11{k}\fi%
\ifnum\titrequatre=12{l}\fi%
\ifnum\titrequatre=13{m}\fi%
}
\long\def\paraquatrel#1{\ifdim\lastskip<\smallskipamount
                \removelastskip\smallskip\fi
                 \vskip0pt plus.01\vsize\penalty-10
                  \vskip0pt
                  plus-.01\vsize{\bigskip
\Indentation{\global\advance\titrequatre by 1
\fontetitrequatre\quad \quad \quad \the\titreunl${\cdot}$\the\titredeuxl${\cdot}$\the\titretrois${\cdot}$\the\titrequatrel.\ }
{\hskip0mm\fontetitrequatre#1}\nobreak\smallskip}}

\ifx\optionkeys\oui
\def\drefun#1{\definexref{¤#1}{{\the\titreun}}{}} 
\def\drefdeux#1{\definexref{¤#1}{{\the\titreun}.{\the\titredeux}}{}}
\def\dreftrois#1{\definexref{¤#1}{{\the\titreun}.{\the\titredeux}.{\the\titretrois}}{}}
\else
\def\drefun#1{\definexref{prg#1}{{\the\titreun}}{}} 
\def\drefdeux#1{\definexref{prg#1}{{\the\titreun}.{\the\titredeux}}{}}
\def\dreftrois#1{\definexref{prg#1}{{\the\titreun}.{\the\titredeux}.{\the\titretrois}}{}}
\fi

%


  \long\def\propdeux#1#2#3#4{%
       \avance{\enonce}
       \leavevmode\edef\temp{#2}%
         \ifx\temp\empty 
          \else
           \definexref{#2}{#1~{\the\titreun.\the\enonce}}{enonces}
            \definexref{s#2}{{\the\titreun.\the\enonce}}{enonces}
             \fi
\smallskip
      \noindent{\bf#1\ {\bf\the\titreun.\the\enonce{#3}.}\enspace}{\sl#4\par}%
      \ifdim\lastskip<\medskipamount \removelastskip\penalty55\par \fi
   }

  \long\def\propun#1#2#3#4{%
      \avance{\enonce}
       \leavevmode\edef\temp{#2}%
        \ifx\temp\empty 
          \else
           \definexref{#2}{#1~{\the\enonce}}{enonces}
            \definexref{{s#2}}{{\the\enonce}}{enonces}
             \fi
   \par 
     \noindent{\bf#1\ {\bf\the\enonce{#3}.}\enspace}{\sl#4\par}%
     \ifdim\lastskip<\medskipamount \removelastskip\penalty55\medskip\fi
  }
  
  \long\def\prop#1#2#3#4{\ifnum\optionparag=1
                          \propdeux{#1}{#2}{\textfont1=\elevenib#3}{#4} \else\propun{#1}{#2}{\textfont1=\elevenib#3}{#4}\fi}

  \long\def\propt#1#2#3{\ifx\tpf\oui \prop{Th\'eo\-r\`eme}{#1}{#2}{#3}\par
                       \else\prop{Theorem}{#1}{#2}{#3}\par\fi}
  \long\def\Propt#1#2{\propt{#1}{}{#2}}
  \long\def\propl#1#2#3{\ifx\tpf\oui\prop{Lem\-me}{#1}{#2}{#3}\par
                         \else\prop{Lemma}{#1}{#2}{#3}\par\fi}
  
  \long\def\propc#1#2#3{\ifx\tpf\oui\prop{Corol\-laire}{#1}{#2}{#3}\par
                         \else\prop{Corollary}{#1}{#2}{#3}\par\fi}

  \long\def\propd#1#2#3{\ifx\tpf\oui\prop{D\'efi\-nition}{#1}{#2}{#3}\par
                       \else\prop{Definition}{#1}{#2}{#3}\par\fi} 
  
  \long\def\proptd#1#2#3{\ifx\tpf\oui\prop{Th\'eor\`eme et d\'efi\-nition}{#1}{#2}{#3}\par
                       \else\prop{Theorem and definition}{#1}{#2}{#3}\par\fi}


  
  \newcount\optionparag\optionparag=1
  
  \long\def\section#1#2{\ifnum\optionparag=1 \paraun{#2} 
                        \else\goodbreak{\fontetitreun
  	                \Indentation{#1.\ }{#2}}\nobreak\smallskip\fi}

  \def\eqconstruct#1{\ifnum\optionparag=1{\the\titreun\hbox{$\cdot$}#1}\else{#1}\fi}

  
  
  \def\numref{oui}  
  
  \newcount\mesref\mesref=0 
  \def\defbib#1{\ifx\numref\oui\global\advance\mesref by 1 \definexref{#1}{{\the
                 \mesref}}{}\else\definexref{#1}{#1}{}\fi}
  \def\bibtem#1{\defbib{#1}\item{\citer{#1}}}
  \def\citer#1{[\ref{#1}]}
  \def\citeplus#1#2{[\ref{#1}; #2]}

  
  \font\seventeenmsa=msam10 at 17pt    
  \font\fourteenmsa=msam10 at 14pt
  \font\twelvemsa=msam10 at 12pt
  \font\tenmsa=msam10                 
  \font\ninemsa=msam10 at 9pt 
  \font\eightmsa=msam10 at 8pt 
  \font\sevenmsa=msam7 
  \font\sixmsa=msam10 at 6pt
  \font\fivemsa=msam5
  \newfam\msafam\textfont\msafam=\tenmsa\scriptfont\msafam=\sevenmsa\scriptscriptfont\msafam=\fivemsa
  
  \font\seventeenbb=msbm10 at 17pt     
  \font\fourteenbb=msbm10 at 14pt
  \font\twelvebb=msbm10 at 12pt
  \font\tenbb=msbm10                   
  \font\ninebb=msbm10 at 9pt 
  \font\eightbb=msbm10 at 8pt 
  \font\sevenbb=msbm7 
  \font\sixbb=msbm10 at 6pt
  \font\fivebb=msbm5 
  \newfam\bbfam\textfont\bbfam=\tenbb\scriptfont\bbfam=\sevenbb\scriptscriptfont\bbfam=\fivebb
  \def\bb{\fam\bbfam\tenbb}%

  \font\seventeenscaln=eusm10 at 17pt   
  \font\twelvescaln=eusm10 at 12pt
  \font\tenscaln=eusm10                
  \font\ninescaln=eusm10 scaled 900
  \font\eightscaln=eusm10 scaled 800
  \font\sevenscaln=eusm10 scaled 700
  \font\sixscaln=eusm10 scaled 600
   
  \newfam\scalnfam\textfont\scalnfam=\tenscaln\scriptfont\scalnfam=\sevenscaln\scriptscriptfont\scalnfam=\sixscaln
  \def\scaln{\fam\scalnfam\tenscaln}%
  \def\scal{\scaln}
  
  \font\tenscalb=eusb10                

  \font\sevenscalb=eusb10 scaled 700

  \newfam\scalbfam\textfont\scalbfam=\tenscalb\scriptfont\scalbfam=\sevenscalb
  %
  
  %
  %
  \font\fourteenrm=cmr12 scaled 1200
  \font\elevenrm=cmr10 at 11pt
  \font\twelverm=cmr12
  \font\ninerm=cmr9
  \font\eightrm=cmr8      
  \font\sevenrm=cmr7
  \font\sixrm=cmr6

  \font\seventeenpcap=cmcsc10 at 17pt
  \font\tenpcap=cmcsc10                        
  \font\ninepcap=cmcsc9
  \font\eightpcap=cmcsc8
  \font\sevenpcap=cmcsc10 scaled 700
  
  \newfam\pcapfam\textfont\pcapfam=\tenpcap\scriptfont\pcapfam=\sevenpcap
  \def\pcap{\fam\pcapfam\tenpcap}
  
  \font\seventeenrm=cmbx12 scaled 1400

  \font\fourteenbf=cmbx10 scaled 1400
  
  \font\twelvebf=cmbx12
  \font\elevenbf=cmbx10 at 11pt
  \font\ninebf=cmbx9  
  \font\eightbf=cmbx8
  \font\sixbf=cmbx6
  
  \font\tengot=eufm10                           
   
  \font\eightgot=eufm10 at 8truept 
  \font\sevengot=eufm7 
  \font\sixgot=eufm10 at 6 truept 
   
  \newfam\gotfam
  \textfont\gotfam=\tengot\scriptfont\gotfam=\sevengot\scriptscriptfont\gotfam=\sixgot
  %

  
  \def\tit{%
  \textfont0=\seventeenrm\scriptfont0=\tenrm\def\rm{\fam0\seventeenrm}%
  \textfont1=\seventeenib\scriptfont1=\twelveib%
  \textfont2=\seventeensy\scriptfont2=\twelvesy\scriptscriptfont2=\ninesy
  \textfont3=\seventeenex
  \textfont\itfam=\seventeenti
  \def\it{\fam\itfam\seventeenti}%
  \textfont\bbfam=\seventeenbb \scriptfont\bbfam=\twelvebb
  \def\bb{\fam\bbfam\seventeenbb}%
  \textfont\msafam=\seventeenmsa\scriptfont\msafam=\twelvemsa
  \textfont\scalnfam=\seventeenscaln
  \def\pcap{\fam\pcapfam\seventeenpcap}
  \normalbaselineskip=25pt\normalbaselines\rm}

  \font\seventeenti=cmbxti10 scaled 1680
  
  \font\fourteenti=cmbxti10 at 14pt
  
  \font\twelveti=cmbxti10 scaled 1200
  \font\eleventi=cmbxti10 at 11pt

  %
  %
  \font\twelveit=cmti12	
  \font\elevenit=cmti10 scaled 1100
  \font\nineit=cmti9
  \font\eightit=cmti8
  \font\sevenit=cmti7

  %
  %
 
 \font\seventeenib=cmmib10 scaled 1680
  \font\fourteenib=cmmib10 scaled 1400
  \font\twelveib=cmmib10 scaled 1200
  \font\elevenib=cmmib10 scaled 1100
  \font\tenib=cmmib10
\font\eightib=cmmib10 scaled 800
  \font\nineib=cmmib10 scaled 900
\font\sevenib=cmmib10 scaled 700
\font\sixib=cmmib10 scaled 600
\font\fiveib=cmmib10 scaled 500

\ifx\ITAN\oui
\else
\innernewfam\cmmibfam
\textfont\cmmibfam=\tenib
\scriptfont\cmmibfam=\sevenib
\scriptscriptfont\cmmibfam=\fiveib
\def\ib{\fam\cmmibfam\tenib}
\fi

  %
  %
  \font\twelvei=cmmi10 scaled 1200
  \font\eleveni=cmmi10 scaled 1100
  \font\ninei=cmmi9
  \font\eighti=cmmi8 
  \font\seveni=cmmi7 			                
  \font\sixi=cmmi6
  
  \font\ninesl=cmsl9                    
  \font\eightsl=cmsl8 
  \font\sevensl=cmsl10 at 7pt

  \font\ninett=cmtt9                    
  \font\eighttt=cmtt8
  \font\seventt=cmtt10 scaled 700

  \font\seventeensy=cmsy10 scaled 1680    
  \font\fourteensy=cmsy10 scaled 1400
  \font\twelvesy=cmsy10 scaled 1176
  
  \font\ninesy=cmsy9                      
  \font\eightsy=cmsy8
  \font\sixsy=cmsy6
  \font\seventeenex=cmex10 at 17pt
  \font\fourteenex=cmex10 at 14pt
  \font\twelveex=cmex10 at 12pt
  \font\nineex=cmex10 at 9pt
  \font\eightex=cmex10 at 8pt
  \font\sevenex=cmex10 at 7pt
  \font\sixex=cmex10 at 6pt
  \font\fiveex=cmex10 at 5pt
  
   
  \font\fourteengp=cmmi10 at 14pt
  
  \font\twelvegp=cmmib10 at 12pt
  \font\elevengp=cmmib10 at 11pt
  \font\tengp=cmmib10                          
  \font\ninegp=cmmib10 at 9pt 
  \font\eightgp=cmmib8 
   
  \font\sixgp=cmmib6


  \def\gponze{\textfont0=\elevenbf\scriptfont0=\eightbf\scriptscriptfont0=\sixbf
           \textfont1=\elevengp\scriptfont1=\eightgp\scriptscriptfont1=\sixgp}
  \def\gpdouze{\textfont0=\twelvebf\scriptfont0=\tenbf\scriptscriptfont0=\ninebf
           \textfont1=\twelvegp\scriptfont1=\tengp\scriptscriptfont1=\ninegp}        
  
 \def\gpquatorze{\textfont0=\fourteenbf\scriptfont0=\twelvebf\scriptscriptfont0=\elevenbf
           \textfont1=\fourteengp\scriptfont1=\twelvegp\scriptscriptfont1=\elevengp}

  
  \expandafter\chardef\csname pre amssym.def at\endcsname=\the\catcode`\@
  \catcode`\@=11
  \def\undefine#1{\let#1\undefined}
  \def\newsymbol#1#2#3#4#5{\let\next@\relax
   \ifnum#2=\@ne\let\next@\msafam@\else
   \ifnum#2=\tw@\let\next@\bbfam@\fi\fi
   \mathchardef#1="#3\next@#4#5}
  \def\mathhexbox@#1#2#3{\relax
   \ifmmode\mathpalette{}{\m@th\mathchar"#1#2#3}%
   \else\leavevmode\hbox{$\m@th\mathchar"#1#2#3$}\fi}
  \def\hexnumber@#1{\ifcase#1 0\or 1\or 2\or 3\or 4\or 5\or 6\or 7\or 8\or
   9\or A\or B\or C\or D\or E\or F\fi}
  
  \def\setboxz@h{\setbox\z@\hbox}
  \def\wdz@{\wd\z@}
  \def\boxz@{\box\z@}
  
  \edef\msafam@{\hexnumber@\msafam}
  \mathchardef\dabar@"0\msafam@39
  
  \edef\bbfam@{\hexnumber@\bbfam}
  \def\widehat#1{\setboxz@h{$\m@th#1$}%
   \ifdim\wdz@>\tw@ em\mathaccent"0\bbfam@5B{#1}%
   \else\mathaccent"0362{#1}\fi}
  \def\widetilde#1{\setboxz@h{$\m@th#1$}%
   \ifdim\wdz@>\tw@ em\mathaccent"0\bbfam@5D{#1}%
   \else\mathaccent"0365{#1}\fi}
  \newsymbol\leqq 1335          
  \newsymbol\leqslant 1336
  \newsymbol\lessgtr 1337       
  \newsymbol\backprime 1038     
  \newsymbol\risingdotseq 133A  
  \newsymbol\fallingdotseq 133B 
  \newsymbol\succcurlyeq 133C   
  \newsymbol\geqq 133D          
  \newsymbol\geqslant 133E
  \newsymbol\nmid 232D
  \newsymbol\nexists 2040
  \newsymbol\smallsetminus 2272
  \newsymbol\varnothing 203F
  
  \catcode`\@=\active

  \catcode`\@=11
  \newcount\typofr\typofr=1
  
  \catcode`\;=\active
  \def;{\ifnum\typofr=1\relax\ifhmode\ifdim\lastskip>\z@\unskip\fi
     \kern.2em\fi\string;\else\string;\fi}
  
  \catcode`\:=\active
  \def:{\ifnum\typofr=1\relax\ifhmode\ifdim\lastskip>\z@\unskip\fi
  \penalty\@M\ \fi\string:\else\string:\fi}
  
  \catcode`\!=\active
  \def!{\ifnum\typofr=1\relax\ifhmode\ifdim\lastskip>\z@\unskip\fi
     \kern.2em\fi\string!\else\string!\fi}
  
  \catcode`\?=\active
  \def?{\ifnum\typofr=1\relax\ifhmode\ifdim\lastskip>\z@\unskip\fi
     \kern.2em\fi\string?\else\string?\fi}

  \def\francais{\typofr=1\def\tpf{oui}}
  \def\anglais{\typofr=2\def\tpf{non}\def\english{oui}}
  \def\oui{oui}
  \francais
  
  \catcode`\@=12
  

\ifx\textures\oui
\def\raye #1|{\leavevmode\setbox1=\hbox{#1}%
\raise .5pt\hbox to \wd1{\xleaders\hbox{\rge{$ \sct / $}%
\kern 1pt}\hfill\kern -1pt }\kern -\wd1 \unhbox1\relax }

\def\barre #1|{\leavevmode\setbox1=\hbox{#1}%
\rlap{\Red\vrule height 2.4pt depth -1.2pt width \wd1}\Black \unhbox1\relax}
\else
\def\raye #1|{\leavevmode\setbox1=\hbox{#1}%
\raise .5pt\hbox to \wd1{\xleaders\hbox{\rge{$ \sct / $}%
\kern 1pt}\hfill\kern -1pt }\kern -\wd1 \unhbox1\relax }

\def\barre #1|{\leavevmode\setbox1=\hbox{#1}%
\rlap{\color{red}\vrule height 2.4pt depth -1.2pt width \wd1}\color{black} \unhbox1\relax}

\fi
  

  
  \def\og{\leavevmode\raise.24ex\hbox{$\scriptscriptstyle\langle\!\langle\>$}}    
  \def\fg{\leavevmode\raise.24ex\hbox{$\scriptscriptstyle\>\rangle\!\rangle$}}    

  \def\d{\,{\rm d}}

  \def\r{{\bb R}}
  \def\CC{{\bb C}}

  \def\F{{\bb F}}
  \def\PP{{\bb P}}

  \def\O{{\scal O}}
  \def\P{{\scaln P}}

  \def\frac#1#2{{#1\over #2}}
  \def\di#1#2{\sct#1\atop{\sct#2}}
  \def\tri#1#2#3{{\sct#1\atop\sct#2}\atop\sct#3}


  \def\¤{\S\thinspace}

  \def\¥{$\bullet$ }
  
  
  \def\e{{\rm e}}
  
  \def\md#1#2{\equiv#1\,({\rm mod\,}#2)}

  \def\epsilon{\varepsilon}

  \def\phi{\varphi}
  \def\theta{\vartheta}
  \def\rho{\varrho}
  \def\dm{{\textstyle{1\over 2}}}

  \def\sct{\scriptstyle}
  \def\pf{\noi{\it Proof. }}
  \def\nid{\ifnum\typofr=1\par\noindent{\it D\'emonstration. }\else\pf\fi}
  \def\noi{\noindent}
  \def\rem{\ifnum\typofr=1\noi{\it Remarque.}\ \else\noi{\it Remark.}\ \fi}
  \def\rems{\ifnum\typofr=1\noi{\it Remarques.}\ \else\noi{\it Remarks.}\ \fi}

  \def\sset{\smallsetminus}

  \def\1{{\bf 1}}
  \def\|{\Vert}

  \def\leq{\leqslant}
  \def\geq{\geqslant}
  \def\wh{\widehat}

  \def\eg{{e.g.}}
  

  \def\fl#1{\left\lfloor #1 \right\rfloor}

  \def\log{\mathop{\rm log}\nolimits}
  
  \def\fs#1#2{{\scriptstyle{#1\over #2}}}
  


\def\abs#1{\left|#1\right|}


  \def\pmb#1{\setbox0=\hbox{#1}%
  \kern-.025em\copy0\kern-\wd0\kern.05em\copy0\kern-\wd0\kern-.025em\raise .0433em\box0 }

  
  \skewchar\eighti='177 \skewchar\sixi='177
  \skewchar\eightsy='60 \skewchar\sixsy='60
  
  \def\eightpoint{%
  \textfont0=\eightrm\scriptfont0=\sixrm\scriptscriptfont0=\fiverm
  \def\rm{\fam0\eightrm}%
  \textfont1=\eighti\scriptfont1=\sixi
  \scriptscriptfont1=\fivei\def\oldstyle{\fam1\seveni}%
  \textfont2=\eightsy\scriptfont2=\sixsy\scriptscriptfont2=\fivesy
  \textfont3=\eightex\scriptfont3=\sixex
  \textfont\itfam=\eightit
  \def\it{\fam\itfam\eightit}%
  \textfont\slfam=\eightsl
  \def\sl{\fam\slfam\eightsl}%
  \textfont\bbfam=\eightbb \scriptfont\bbfam=\sixbb\scriptscriptfont\bbfam=\fivebb
  \def\bb{\fam\bbfam\eightbb}%
  \textfont\msafam=\eightmsa\scriptfont\msafam=\sixmsa
  \textfont\scalnfam=\eightscaln
  \def\scaln{\fam\scalnfam\eightscaln}
  \textfont\ttfam=\eighttt
  \def\tt{\fam\ttfam\eighttt}%
\textfont\gotfam=\eightgot
  \textfont\bffam=\eightbf\scriptfont\bffam=\sixbf\scriptscriptfont\bffam=\fivebf
  \def\bf{\fam\bffam\eightbf}%
  \ifx\ITAN\oui\else\textfont\cmmibfam=\eightib
       \scriptfont\cmmibfam=\sixib
        \scriptscriptfont\cmmibfam=\fiveib
         \def\ib{\fam\cmmibfam\eightib}
   \fi
  \textfont\pcapfam=\eightpcap
  \def\pcap{\fam\pcapfam\eightpcap}
  \abovedisplayskip=2pt plus2pt minus 2pt
  \belowdisplayskip=2pt plus1pt minus 2pt
  \abovedisplayshortskip= 1pt plus 2pt minus 1pt
  \belowdisplayshortskip= 1pt plus 2pt minus 1pt
  \smallskipamount=2pt plus 1pt minus 2pt
  \medskipamount=3pt plus 2pt minus 2pt
  \bigskipamount=7pt plus 3pt minus 3pt
  \setbox\strutbox=\hbox{\vrule height 5pt depth 2pt width 0pt}%
  \normalbaselineskip=9pt\normalbaselines\rm}

  \def\({\left(}
  \def\){\right)}
  
  \def\footnoterule{\kern -2pt\hrule width 7truecm\kern 2.4pt}
  
  \def\xnotedef#1{\definexref{#1}{\noexpand\number\footnotenumber}{Note}}%

  
  
  \def\ninepoint{%
  \textfont0=\ninerm\scriptfont0=\sixrm\scriptscriptfont0=\fiverm
  \def\rm{\fam0\ninerm}%
  \textfont1=\ninei\scriptfont1=\sixi
  \scriptscriptfont1=\fivei\def\oldstyle{\fam1\ninei}%
  \textfont2=\ninesy\scriptfont2=\sixsy\scriptscriptfont2=\fivesy
  \textfont3=\nineex\scriptfont3=\sixex
  \textfont\itfam=\nineit
  \def\it{\fam\itfam\nineit}%
  \textfont\slfam=\ninesl
  \def\sl{\fam\slfam\ninesl}%
  \textfont\bbfam=\ninebb\scriptfont\bbfam=\sixbb\scriptscriptfont\bbfam=\fivebb
  \def\bb{\fam\bbfam\ninebb}%
  \textfont\msafam=\ninemsa\scriptfont\msafam=\sixmsa\scriptscriptfont\msafam=\fivemsa
  \textfont\scalnfam=\ninescaln
  \def\scaln{\fam\scalnfam\ninescaln}
  \textfont\ttfam=\ninett
  \def\tt{\fam\ttfam\ninett}%
  \textfont\bffam=\ninebf\scriptfont\bffam=\sixbf\scriptscriptfont\bffam=\fivebf
  \def\bf{\fam\bffam\ninebf}%
  \abovedisplayskip=3pt plus2pt minus 2pt
  \belowdisplayskip=3pt plus1pt minus 2pt
  \abovedisplayshortskip= 2pt plus 2pt minus 1pt
  \belowdisplayshortskip= 2pt plus 2pt minus 1pt
  \smallskipamount=2pt plus 1pt minus 2pt
  \medskipamount=3pt plus 2pt minus 2pt
  \bigskipamount=7pt plus 3pt minus 3pt
  \setbox\strutbox=\hbox{\vrule height 5pt depth 2pt width 0pt}%
  \normalbaselineskip=11pt plus.3pt minus.2pt\normalbaselines\rm}

  \def\sevenpoint{%
  \textfont0=\sevenrm\scriptfont0=\sixrm\scriptscriptfont0=\fiverm
  \def\rm{\fam0\sevenrm}%
  \textfont1=\seveni\scriptfont1=\sixi
  \scriptscriptfont1=\fivei\def\oldstyle{\fam1\seveni}%
  \textfont2=\sevensy\scriptfont2=\sixsy\scriptscriptfont2=\fivesy
  \textfont3=\sevenex\scriptfont3=\fiveex
  \textfont\itfam=\sevenit
  \def\it{\fam\itfam\sevenit}%
  \textfont\slfam=\sevensl
  \def\sl{\fam\slfam\sevensl}%
  \textfont\bbfam=\sevenbb \scriptfont\bbfam=\sixbb\scriptscriptfont\bbfam=\fivebb
  \def\bb{\fam\bbfam\sevenbb}%
  \textfont\msafam=\sevenmsa\scriptfont\msafam=\sixmsa
  \textfont\scalnfam=\sevenscaln
  \def\scaln{\fam\scalnfam\sevenscaln}
  \textfont\bffam=\sevenbf\scriptfont\bffam=\sixbf\scriptscriptfont\bffam=\fivebf
  \def\bf{\fam\bffam\sevenbf}%
  \textfont\ttfam=\seventt
  \abovedisplayskip=2pt plus2pt minus 2pt
  \belowdisplayskip=2pt plus1pt minus 2pt
  \abovedisplayshortskip= 1pt plus 2pt minus 1pt
  \belowdisplayshortskip= 1pt plus 2pt minus 1pt
  \smallskipamount=2pt plus 1pt minus 2pt
  \medskipamount=3pt plus 2pt minus 2pt
  \bigskipamount=7pt plus 3pt minus 3pt
  \setbox\strutbox=\hbox{\vrule height 5pt depth 2pt width 0pt}%
  \normalbaselineskip=9pt\normalbaselines\rm}

 \def\onzepts{%
 \textfont0=\elevenrm\scriptfont0=\ninerm
 \textfont1=\eleveni\scriptfont1=\ninei
}

\def\douzepts{%
  \textfont0=\twelverm\scriptfont0=\tenrm\def\rm{\fam0\twelverm}%
  \textfont1=\twelvei\scriptfont1=\teni%
  \textfont2=\twelvesy\scriptfont2=\tensy\scriptscriptfont2=\eightsy
  \textfont3=\twelveex
  \textfont\itfam=\twelveti
  \def\it{\fam\itfam\twelveti}%
  \textfont\bffam=\twelvebf\scriptfont\bffam=\tenbf\scriptscriptfont\bffam=\eightbf
  \def\bf{\fam\bffam\twelvebf}%
  \textfont\bbfam=\twelvebb \scriptfont\bbfam=\tenbb
  \def\bb{\fam\bbfam\twelvebb}%
  \textfont\msafam=\twelvemsa\scriptfont\msafam=\tenmsa
  \textfont\scalnfam=\twelvescaln
  \normalbaselineskip=15pt\normalbaselines\rm}

\def\quatorzepts{%
  \textfont0=\fourteenrm\scriptfont0=\twelverm\def\rm{\fam0\fourteenrm}%
  \textfont1=\fourteenib\scriptfont1=\twelveib%
  \textfont2=\fourteensy\scriptfont2=\twelvesy\scriptscriptfont2=\tensy
  \textfont3=\fourteenex
  \textfont\itfam=\fourteenti
  \def\it{\fam\itfam\fourteenti}%
  \textfont\bffam=\fourteenbf\scriptfont\bffam=\twelvebf\scriptscriptfont\bffam=\tenbf
  \def\bf{\fam\bffam\fourteenbf}%
  \textfont\bbfam=\fourteenbb \scriptfont\bbfam=\twelvebb
  \def\bb{\fam\bbfam\fourteenbb}%
  \textfont\msafam=\fourteenmsa\scriptfont\msafam=\twelvemsa
  \textfont\scalnfam=\twelvescaln
  \normalbaselineskip=18pt\normalbaselines\rm}


\def\AA{{\it Acta Arith.}}

\def\picture #1 by #2 (#3){\leavevmode\vbox to #2{
     \hrule width #1 height 0pt depth 0pt
      \vfill
       \special{picture #3}}}

\def\illustration #1 by #2 (#3) scaled #4{\dimen1=#2
  \divide\dimen1 by 1000
  \multiply\dimen1 by #4
  \vtop to \dimen1{\dimen1=#1
  \divide\dimen1 by 1000
  \multiply\dimen1 by #4
  \hsize=\dimen1\vss
  \noindent\special{illustration #3 scaled #4}}}

\ifx\couleurs\oui

\fi

%% file: option_keys.tex
\ifx\optionkeymacros\undefined\else \fi

\catcode`\Œ=\active\defŒ{{\aa}}       
\catcode`\º=\active\defº{\int}        
\catcode`\=\active\def{\c c}        
\catcode`\¶=\active\def¶{\partial}    
\catcode`\Ä=\active\defÄ{\oint}       
\catcode`\Æ=\active\defÆ{\triangle}   
\catcode`\Â=\active\defÂ{\neg}        
\catcode`\µ=\active\defµ{\mu}         
\catcode`\¿=\active\def¿{{\o}}        
\catcode`\¹=\active\def¹{\pi}         
\catcode`\Ï=\active\defÏ{{\oe}}       
\catcode`\§=\active\def§{{\ss}}       
\catcode`\ =\active\def {\dagger}     
\catcode`\Ã=\active\defÃ{\sqrt}       
\catcode`\·=\active\def·{\Sigma}      
\catcode`\Å=\active\defÅ{\approx}     
\catcode`\½=\active\def½{\Omega}      
\catcode`\£=\active\def£{{\it\$}}     
\catcode`\°=\active\def°{\infty}      
\catcode`\¤=\active\def¤{{\S}}        
\catcode`\¦=\active\def¦{{\P}}        
\catcode`\¥=\active\def¥{\bullet}     
\catcode`\»=\active\def»{\leavevmode\raise.585ex\hbox{\b a}}      
\catcode`\¼=\active\def¼{\leavevmode\raise.6ex\hbox{\b o}}        
\catcode`\­=\active\def­{\not=}       
\catcode`\²=\active\def²{\leq}        
\catcode`\³=\active\def³{\geq}        
\catcode`\Ö=\active\defÖ{\div}        
\catcode`\É=\active\defÉ{{\dots}}     
\catcode`\¾=\active\def¾{{\ae}}       
\catcode`\Ç=\active\defÇ{\og}         
\catcode`\Ò=\active\defÒ{``}          
\catcode`\Á=\active\defÁ{!`}          
\catcode`\¢=\active\def¢{\rlap/c}     
\catcode`\Ô=\active\defÔ{`}           
\catcode`\Õ=\active\defÕ{'}           


\catcode`\=\active\def{{\AA}}       
\catcode`\'=\active\def'{\c C}        
\catcode`\¯=\active\def¯{{\O}}        
\catcode`\¸=\active\def¸{\Pi}         
\catcode`\Î=\active\defÎ{{\OE}}       
\catcode`\®=\active\def®{{\AE}}       
\catcode`\×=\active\def×{\diamond}    
\catcode`\¡=\active\def¡{\accent'27}  
\catcode`\Ó=\active\defÓ{''}          
\catcode`\±=\active\def±{\pm}         
\catcode`\È=\active\defÈ{\fg}         
\catcode`\À=\active\defÀ{?`}          
\catcode`\Ð=\active\defÐ{--}          
\catcode`\Ñ=\active\defÑ{---}         


\catcode`\Š=\active\defŠ{\"a}        
\catcode`\'=\active\def'{\"e}        
\catcode`\•=\active\def•{\"{\i}}     
\catcode`\š=\active\defš{\"o}        
\catcode`\Ÿ=\active\defŸ{\"u}        
\catcode`\Ø=\active\defØ{\"y}        
\catcode`\å=\active\defå{\^A}        
\catcode`\€=\active\def€{\"A}        
\catcode`\…=\active\def…{\"O}        
\catcode`\†=\active\def†{\"U}        
\catcode`\‡=\active\def‡{\'a}        
\catcode`\Ž=\active\defŽ{\'e}        
\catcode`\'=\active\def'{\'{\i}}     
\catcode`\—=\active\def—{\'o}        
\catcode`\œ=\active\defœ{\'u}        
\catcode`\ƒ=\active\defƒ{\'E}        
\catcode`\æ=\active\defæ{\^E}        
\catcode`\é=\active\defé{\`E}        %
\catcode`\ˆ=\active\defˆ{\`a}        
\catcode`\=\active\def{\`e}        
\catcode`\"=\active\def"{\`{\i}}     
\catcode`\˜=\active\def˜{\`o}        
\catcode`\=\active\def{\`u}        
\catcode`\Ë=\active\defË{\`A}        
\catcode`\‹=\active\def‹{\~a}        
\catcode`\–=\active\def–{\~n}        
\catcode`\›=\active\def›{\~o}        
\catcode`\Ì=\active\defÌ{\~A}        
\catcode`\"=\active\def"{\~N}        
\catcode`\Í=\active\defÍ{\~O}        
\catcode`\‰=\active\def‰{\^a}        
\catcode`\=\active\def{\^e}        
\catcode`\"=\active\def"{\^{\i}}     
\catcode`\™=\active\def™{\^o}        
\catcode`\ž=\active\defž{\^u}        

\let\optionkeymacros\null